\documentclass[a4paper,12pt,reqno]{amsart}
\usepackage{latexsym}
\usepackage{amssymb}
\usepackage{amsmath}
\def\spann{\mathrm{span}}
\def\Id{\mbox{Id}}
\def\id{\mbox{id}}

\def\vol{\mbox{vol}}
\def\End{\mbox{End}}
\def\Re{\mbox{Re}}
\def\tr{\mbox{tr}}
\def\Sym{\mbox{Sym}}

\def\a{\alpha}

\def\beq{\begin{equation}}
\def\eeq{\end{equation}}
\def\bea{\begin{eqnarray*}}
\def\eea{\end{eqnarray*}}

\def\i{\lrcorner\,}

\def\la{\langle}

\def\n{\nabla}

\def\o{\omega}

\def\SU{{\rm SU}}

\def\Id{{\rm Id}}

\def\r{\rightline{$\Box$}\ss}
\def\ra{\rangle}
\def\ss{\smallskip}
\def\t{\tilde}

\def\L{{\mathcal L}}
\def\K{{\mathcal K}}
\def\O{\Omega}

\def\Ric{{\rm Ric}}

\def\CM{\ensuremath{\mathbb C}}

\def\s{\sigma}

\def\.{\cdot}

\def\J{\widehat{J}}

\newcommand{\firstline}[2]{\makebox[#1][l]{$\displaystyle{#2}$}&&}
\newtheorem{ede}{Definition}[section]

\newtheorem{epr}[ede]{Proposition}

\newtheorem{ath}[ede]{Theorem}

\newtheorem{elem}[ede]{Lemma}

\newtheorem{ecor}[ede]{Corollary}
 
\title{Unit Killing Vector Fields on Nearly K{\"a}hler Manifolds}

\author{Andrei Moroianu, Paul--Andi Nagy and Uwe Semmelmann}

\address{Andrei Moroianu \\ CMAT\\ {\'E}cole Polytechnique \\ UMR 7640 du CNRS
\\ 91128 Palaiseau \\ France}
\email{am@math.polytechnique.fr}

\address{Paul--Andi Nagy \\
Institut f{\"u}r Mathematik \\
Humboldt Universit{\"a}t zu Berlin \\ 
Sitz : Rudower Chaussee 25, D-12489, Berlin, Germany}
\email{nagy@mathematik.hu-berlin.de}

\address{Uwe Semmelmann\\ Fachbereich Mathematik, Universit{\"a}t
Hamburg\\ Bundesstr. 55,  D-20146 Hamburg, Germany}
\email{Uwe.Semmelmann@math.uni-hamburg.de}

\begin{document}

\begin{abstract} We study $6$-dimensional nearly K{\"a}hler manifolds
  admitting a Killing vector field of unit length. In the compact
case it is shown that up to a finite cover there is only one geometry possible, that of the
$3$--symmetric space $S^{3} \times S^3$.

\end{abstract}

\maketitle

\section{Introduction}
Nearly K{\"a}hler geometry (shortly NK in what follows) naturally arises
as one of the sixteen classes 
of almost Hermitian manifolds appearing in the celebrated
Gray--Hervella classification \cite{Hervella}. These 
manifolds were studied intensively in the seventies by A. Gray
\cite{Gray}. His initial motivation was 
inspired by the concept of weak holonomy \cite{Gray}, but very
recently it turned out that this concept, as defined by 
Gray, does not produce any new geometric structure (see
\cite{Alexandrov}) other than those coming from a Riemannian holonomy
reduction. 
One of the most important
properties of NK manifolds is that their canonical Hermitian
connection has totally skew--symmetric, parallel 
torsion \cite{Kiri}. From this point of view, they naturally fit into the
setup proposed in \cite{Cleyton-Swann} towards a weakening of the 
notion of Riemannian holonomy. The same property suggests that
NK manifolds might be objects of interest in string theory
\cite{friedrich2}.

The structure theory of compact NK manifolds, as developed in
\cite{Nagy} reduces their study to positive quaternionic K{\"a}hler
manifolds and nearly--K{\"a}hler manifolds of dimension $6$.
The last class of manifolds falls in the area of special metrics with
very rigid -- though not yet fully understood -- properties.

Indeed, it is known since a long time that in $6$ dimensions, a
NK metric which is not K{\"a}hler has to be Einstein of positive scalar
curvature. Moreover, such a structure is characterized by the
existence of some (at least locally defined) real Killing spinor
\cite{Grunewald}.
Combining these properties with the fact that the first Chern class
(form) vanishes \cite{Gray}, one observes that non--K{\"a}hler,
nearly K{\"a}hler
six--dimensional manifolds solve most of the type II string equations
\cite{friedrich2}. Despite of all these interesting features, very
little is known about these manifolds. In particular, apart from the
$3$--symmetric spaces
$$ S^6, S^3 \times S^3, \mathbb{C}P^3, F(1,2) $$
no compact example is available and moreover these is are the only
compact homogeneous examples \cite{Butruille}. 

In a recent article \cite{Hitchin}, Hitchin shows that
nearly parallel  $G_2$--structures (\cite{fkms} for an account) 
and NK manifolds of $6$ dimensions have the same variational
origins. On the other hand, 
many examples of nearly parallel $G_2$--structures are available 
since any $7$--dimensional, $3$--Sasakian manifold carries such a
structure \cite{fkms, GaSa} and a profusion 
of compact examples of the latter were produced in \cite{Galicki}. 
Since one property of $3$--Sasakian manifolds is to admit
unit Killing vector fields, one might ask whether this can happen in the
NK setting. 

In the present paper we study $6$--dimensional non--K{\"a}hler,
nearly K{\"a}hler manifolds which globally admit a Killing vector field
$\xi$ of
constant length. After recalling 
some elementary features of nearly K{\"a}hler geometry in
Section 2, we show in Section 3 that any Killing vector field 
of unit length induces a transversal almost hyper--Hermitian structure
on the manifold.
The almost hyper--Hermitian structure is preserved by the Killing
vector field $\xi$ but not by 
$J \xi$ (here $J$ denotes the almost complex structure of the
nearly K{\"a}hler structure). We measure this in the fourth section by 
computing the Lie derivatives of the various geometrically significant
tensors in the direction of $J \xi$. This technical part is used 
in Section 5 to perform a double reduction of the $6$--dimensional
nearly K{\"a}hler manifold. The resulting $4$--dimensional manifold 
is in fact a K{\"a}hler--Einstein surface of positive scalar curvature
admitting a orthogonal almost--K{\"a}hler structure inducing the opposite
orientation. 
The geometry of the situation is completely understood in terms of
this data. Moreover, if the nearly K{\"a}hler manifold is compact, a
Sekigawa--type argument from \cite{adm} 
shows that the almost--K{\"a}hler structure is
actually integrable, allowing us to prove the main result of this paper. 
\begin{ath}\label{main}
Let $(M^6,g,J)$ be a complete nearly K{\"a}hler manifold. If $g$ admits a unit
Killing vector field, then up to a finite cover $(M^6,g,J)$ is
isometric to $S^3 \times S^3$ endowed with its canonical NK structure.
\end{ath}

\section{Nearly K{\"a}hler Manifolds}

An almost Hermitian manifold $(M, g, J)$ is called
{\it nearly K{\"a}hler} if $(\nabla_X J) X = 0$ is satisfied
for all vector fields $X$. In other words, the covariant derivative of
$J$ (viewed as a ($3,0$)--tensor via the metric $g$) 
is skew--symmetric in all three arguments, not only in the last two,
as it is the case for general almost Hermitian structures. 
This is equivalent to 
$d \Omega = 3 \nabla \Omega$, where $\Omega$ is the fundamental 
$2$--form, {\em i.e.} $\Omega(X, Y) := g(J X, Y)$. The following 
lemma summarizes some of the known identities for nearly K{\"a}hler manifolds.

\begin{elem}\label{nk}
Let $(M, g, J)$ be a nearly K{\"a}hler manifold. Then
\begin{enumerate}
\item \quad $(\nabla_X J) Y + (\nabla_Y J) X = 0 $
\item \quad $(\nabla_{JX} J) Y = (\nabla_{X} J) JY $
\item \quad 
$
J((\nabla_X J) Y) = - (\nabla_X J) JY = - (\nabla_{JX} J) Y
$
\item \quad
$ g(\nabla_X Y, X) = g(\nabla_X JY, JX) $
\item \quad
$
2 g((\nabla^2_{W, X} J) Y, Z)
 = - \sigma_{X,Y,Z}   g((\nabla_W J) X,  (\nabla_Y J) JZ) 
$
\end{enumerate}
where $ \sigma_{X,Y,Z}  $ denotes the cyclic sum over 
the vector fields $ X, Y, Z$.
\end{elem}

A nearly K{\"a}hler manifold is called
to be of {\it constant type} $\alpha$ if
$$
\|(\nabla_X J)(Y)\|^2 = \alpha \{\|X\|^2 \|Y\|^2 - g(X, Y)^2
 - g(J X, Y)^2\}
$$
holds for any vector fields $X, Y$. A. Gray proved that a
nearly K{\"a}hler manifold of positive constant type is necessarily
$6$--dimensional (cf.~\cite{Gray}). Moreover he showed:
\begin{epr}
Let $(M, g, J)$ be a $6$--dimensional nearly K{\"a}hler, non--K{\"a}hler
manifold, then
\begin{enumerate}
\item\quad $M$ is of constant type $\alpha > 0$.
\item\quad $c_1(M) = 0$ and in particular $M$ is a spin manifold.
\item\quad $(M, g)$ is Einstein and $\Ric=5 \alpha \Id=5 \Ric^*$
\end{enumerate}
\end{epr}
Here the $\ast$--Ricci curvature $\Ric^*$ is defined as
$\Ric^*(X, Y) = \tr(Z\mapsto R(X, JZ) JY))$. From this it easily
follows that $\Ric^*(X, Y) = {\mathcal R}(\Omega)(X, JY)$, where $\mathcal R$
denotes the curvature operator on $2$--forms.

\begin{elem}
Let $(M^6, g J)$ be a nearly K{\"a}hler manifold of constant type $\alpha$, then
\bea
g((\nabla_U J) X, (\nabla_Y J) Z)
&=&
\alpha \{
g(U, Y) g(X, Z) - g(U, Z) g(X, Y)\\
&&\phantom{xxxxxxxxx} - 
g(U, JY) g(X, JZ) + g(U, JZ) g(X, JY) \}
\eea
\end{elem}
\begin{ecor}\label{nabla}
Let $(M^6, g, J)$ be a nearly K{\"a}hler manifold of constant type 
$\alpha=1$, then
\bea
(\nabla_X J)\circ (\nabla_X J)  Y
&=&
- |X|^2 Y\qquad \qquad\mbox{for} \quad Y \perp X, JX\\[1.5ex]
\nabla^* \nabla \Omega
&=&
\quad 4 \Omega
\eea
\end{ecor}

\begin{elem}\label{ortho}
Let $X$ and $Y$ be any vector fields on $M$, then the vector field
$(\nabla_XJ) Y$ is orthogonal to $X, JX, Y, $ and $JY$.
\end{elem}

This lemma allows us to use adapted frames $\{e_i\}$ which
are especially convenient for local calculations. Let $e_1$ and
$e_3$ be any two orthogonal vectors and define:
$$
e_2:=Je_1,\quad e_4:=Je_3,\qquad
e_5:=(\nabla_{e_1} J) e_3,\qquad
e_6:=Je_5
$$
\begin{elem}
With respect to an adapted frame $\{e_i\}$ one has
\bea
\nabla J &=& \quad e_1 \wedge e_3 \wedge e_5 - 
e_1 \wedge e_4 \wedge e_6 - 
e_2 \wedge e_3 \wedge e_6 - 
e_2 \wedge e_4 \wedge e_5\\
\ast (\nabla J) &=& -  e_2 \wedge e_4 \wedge e_6 + 
e_2 \wedge e_3 \wedge e_5 + 
e_1 \wedge e_4 \wedge e_5 + 
e_1 \wedge e_3 \wedge e_6\\
\Omega &=& \quad e_1 \wedge e_2 + e_3 \wedge e_4 + e_5 \wedge e_6
\eea
\end{elem}
\begin{ecor}\label{elementary}
Let $\Omega$ be the fundamental $2$--form and $X$ an arbitrary vector
field, then 
\begin{enumerate}
\item\quad
$
X \i \Omega = J X^{\flat},\quad
X \i \ast \Omega = JX^{\flat} \wedge \Omega,\quad
X \i d \Omega = JX \i \ast d \Omega
$
\item\quad
$
| \Omega |^2 = 3,\quad
\ast \Omega = \frac12 \Omega \wedge \Omega,\quad
\vol = e_1 \wedge\ldots\wedge e_6 = \frac16 \Omega^3,\quad
\Omega \wedge d \Omega = 0
$
\item\quad
$
\ast X^{\flat} = \frac12 JX^{\flat} \wedge \Omega \wedge \Omega
,$

where $X^{\flat}$ denotes the $1$--form which is metric dual to $X$.
\end{enumerate}

\end{ecor}

\begin{epr}
Let $(M^6, g, J)$ be a nearly K{\"a}hler, non--K{\"a}hler manifold
with fundamental $2$--form $\Omega$ which is of constant type 
$\alpha = 1$, then
$$
\Delta \Omega = 12 \Omega
$$
\end{epr}
{\sc Proof.}
To show that $\Omega$ is an eigenform of the Laplace operator
we will use the Weitzenb{\"o}ck formula on $2$--forms, {\em i.e.} 
$\Delta  = \nabla^*\nabla + \frac{s}{3} \Id-2 {\mathcal R}$.
From Corollary \ref{nabla} we know that 
$\nabla^*\nabla \Omega = 4 \Omega$. 
Since we assume $M$ to be of constant type $1$ the scalar curvature
is $s=30$ and the $\ast$--Ricci curvature $\Ric^*$ coincides with
the Riemannian metric $g$. Hence,
$ {\mathcal R}(\Omega)(X, Y) = - Ric^*(X, J Y) = \Omega(X, Y)$. 
Substituting this into the Weitzenb{\"o}ck formula yields $\Delta\Omega=12\Omega$.

\r

\begin{ecor}\label{djxi}
If $\xi$ is any vector field satisfying $L_\xi( \ast  d  \Omega)=0$
then
\beq\label{c29}
d (J\xi \i d\Omega)=- 12 J\xi^{\flat} \wedge \Omega
\eeq
\end{ecor}
{\sc Proof.}
We use Corollary \ref{elementary} and the relation 
$L_\xi = d \circ\xi \i + \xi \i \circ d$ to obtain
$$
d(J\xi \i d\Omega)
 = 
- d(\xi \i \ast d \Omega)
 = 
\xi \i (d \ast d \Omega)
 = 
- \xi \i (\ast d^* d \Omega)
 = 
- \xi \i (\ast \Delta \Omega)
$$
Note that $d^* \Omega = 0$. Again applying Corollary \ref{elementary}
together with (\ref{c29}) we get
$$
d(J\xi \i d\Omega)
 = 
- 12 \xi \i (\ast \Omega)
 = 
-6 \xi \i (\Omega \wedge \Omega)
 = 
- 12 (\xi \i \Omega) \wedge \Omega
 = 
- 12 J\xi^{\flat}  \wedge \Omega
$$

\r

The structure group of a nearly K{\"a}hler manifold $(M^6, g, J)$
reduces to $\SU(3)$ which implies the decomposition
$\Lambda^2(TM)=\Lambda^{inv}\oplus\Lambda^{anti}$ with
\bea
\alpha\in\Lambda^{inv} 
&\Leftrightarrow &
\alpha (X, Y) = \quad\alpha (JX, JY)\\
\alpha\in\Lambda^{anti}
&\Leftrightarrow &
\alpha (X, Y) = - \alpha (JX, JY)
\eea
We will denote the projection of a $2$--form $\alpha$ onto $\Lambda^{inv}$
by $\alpha^{(1,1)}$ and the projection onto $\Lambda^{anti}$ by
$\alpha^{(2,0)}$. This is motivated by the isomorphisms
$$
\Lambda^{inv}\otimes\CM \cong \Lambda^{(1,1)}(TM) \cong {\rm u}(3),
\qquad
\Lambda^{anti}\otimes\CM \cong 
\Lambda^{(2,0)}(TM) \oplus \Lambda^{(0,2)}(TM)
$$
\begin{elem}
The decomposition $\Lambda^2(TM)=\Lambda^{inv}\oplus\Lambda^{anti}$
is orthogonal and the projections of a $2$--form $\alpha$ onto the
two components are given by
\bea
\alpha^{(1,1)} (X, Y)
&=&
\frac12 (\alpha (X, Y) + \alpha (JX, JY))
 = 
\frac12 \Re  (\alpha([X+iJX], [Y-iJY]))\\
\alpha^{(2,0)} (X, Y)
&=&
\frac12 (\alpha (X, Y) - \alpha (JX, JY))
 = 
\frac12 \Re  (\alpha([X+iJX], [Y+iJY]))
\eea
\end{elem}

Under the isomorphism $\Lambda^2(TM)\stackrel{\sim}{\rightarrow} \End_0(TM)$
every $2$--form $\alpha$ corresponds to a skew--symmetric endomorphism
$A$ which is defined by the equation $\alpha(X, Y) = g(AX, Y)$.
Note that $|A|^2 = 2 |\alpha|^2$, where $|\cdot |$ is the norm
induced from the Riemannian metric on the endomorphisms and on the
$2$--forms.

\begin{elem}\label{invanti}
Let $A^{(2,0)}$ resp. $A^{(1,1)}$ be the endomorphisms corresponding
to the components $\alpha^{(2,0)}$ resp. $\alpha^{(1,1)}$, then
\begin{enumerate}
\item \quad
$
A^{(2,0)} = \frac12 (A + J A J),\qquad
A^{(1,1)} = \frac12 (A - J A J)
$
\item \quad
$
J\circ A^{(1,1)} = A^{(1,1)}\circ J,\qquad
J\circ A^{(2,0)} = - A^{(2,0)}\circ J
$
\end{enumerate}
\end{elem}


\section{The Transversal Complex Structures}

In this section we consider  $6$--dimensional compact non--K{\"a}hler
nearly K{\"a}hler manifolds $(M, g, J)$ of constant type $\alpha
\equiv 1$, {\it i.e.} with scalar curvature normalized by
$s \equiv 30$. We will assume from now on that $ (M, g) $
is not isometric to the sphere with its standard metric.

\begin{epr}
Let $\xi$ be a Killing vector field then the Lie derivative of
the almost complex structure $J$ with respect to $\xi$ vanishes
$$
L_\xi J = 0.
$$
\end{epr}
{\sc Proof.}
Nearly K{\"a}hler structures $J$ on $(M, g)$ are in one--to--one correspondence
with Killing spinors of unit norm on $M$ (cf.~\cite{Grunewald}). 
However, if  $ (M, g) $ is not isometric to the standard  sphere,
the (real) space of Killing spinors is $1$--dimensional, so there
exist exactly $2$ nearly K{\"a}hler structures compatible with $g$: $J$
and $-J$. This shows that the identity component of the isometry group
of $M$ preserves $J$, so in particular $ L_\xi J= 0 $ for every
Killing vector field $\xi$.

\r

Since a Killing vector field $ \xi $ satisfies by definition
$ L_\xi g= 0 $ we obtain that the Lie derivative $ L_\xi $
of all natural tensors constructed out of $g$ and $J$ vanishes.
In particular, we have

\begin{ecor}\label{enc} 
If $\xi$ is a Killing vector field on the
  nearly K{\"a}hler manifold $(M,g,J)$ with fundamental $2$--form
  $\Omega$, then
$$
L_\xi (\Omega)  =  0,\qquad
L_\xi (d \Omega)  =  0,\qquad
L_\xi (\ast d \Omega)  =  0.\qquad
$$
\end{ecor}

In order to simplify the notations we denote by $\zeta:=\xi^{\flat}$
the dual 1--form to $\xi$. 

\begin{ecor}\label{dxi}
$\qquad
d(J \zeta) = - \xi \i d\Omega
$
\end{ecor}
{\sc Proof.}
From $L_\xi \Omega = 0$ follows: $ - \xi \i d\Omega
=d (\xi \i \Omega)=d(J \zeta)$.

\r

\bigskip

\begin{elem}\label{dxi20}
Let $\xi$ be a Killing vector field with metric dual $\zeta$ and let
$d\zeta = d\zeta^{(1,1)} + d\zeta^{(2,0)}$ be the type decomposition of
$d\zeta$. Then the endomorphisms corresponding to $ d\zeta^{(2,0)} $ and
$ d\zeta^{(1,1)} $ are $ - \nabla_{J\xi}J $ and
$ 2\nabla_{\cdot}\xi +  \nabla_{J\xi}J$ respectively.
\end{elem}
{\sc Proof.}
The equation $ L_\xi J = 0$ applied to a vector field 
$X$ yields $[\xi, JY]=J [\xi, Y]$, which can be written as
$$
\nabla_\xi JY - \nabla_{JY} \xi = J \nabla_\xi Y - J \nabla_{Y} \xi.
$$
From this equation we obtain
\beq\label{rt}
(\nabla_\xi J) Y =  \nabla_{JY}\xi - J \nabla_{Y} \xi.
\eeq
Let $A$ ($=2\nabla_{\cdot}\xi$) denote the skew--symmetric endomorphism
corresponding to the $2$--form $d\zeta$.
Then (\ref{rt}) can be written 
$$
\nabla_\xi J = \frac12  (A\circ J -  J\circ A)
$$
From Lemma \ref{nk} we  get
$  \nabla_{J\xi} J = (\nabla_{\xi} J)\circ J$, whence
$$
\nabla_{J\xi} J
 = 
\frac12  (A\circ J -  J\circ A)\circ J
 = 
-  \frac12 (A +  J\circ A\circ J)
 = - A^{(2,0)}
$$
For the corresponding $2$--form we obtain: 
$d \zeta^{(2,0)} = -\nabla_{J\xi} \Omega$. Finally
we compute $ A^{(1,1)} $ as
$$
A^{(1,1)}
  =  
A   -  A^{(2,0)}
  =  
2 \nabla_{\cdot}\xi   +  \nabla_{J\xi} J \ .
$$

\r

From now on we will mainly be interested in compact nearly K{\"a}hler
manifolds admitting a Killing vector field $\xi$ of constant length
(normalized to $1$). We start by collecting several elementary properties

\begin{elem}\label{l35}
The following relations hold:
\begin{enumerate}
\item \quad
$\nabla_{\xi} \xi = 
\nabla_{J\xi} \xi = 
\nabla_{\xi} J\xi = 
\nabla_{J\xi} J\xi = 0,\qquad [\xi, J\xi] = 0
$
\item \quad
$\xi  \i  d\zeta = J\xi \i  d\zeta = 0 \ . $
\end{enumerate}
In particular, the distribution $ V:=\spann\{\xi, J\xi\} $ is integrable.
\end{elem}

We now define two endomorphisms which
turn out to be complex structures on the orthogonal complement
$  H:= \{\xi, J\xi\}^\perp$:
\beq\label{ik}
I  :=   \nabla_\xi J,
\qquad\mbox{and}\qquad
K  :=  \nabla_{J\xi} J \ .
\eeq

Note that $ -K $ is the endomorphism corresponding to
$ d\zeta^{(2,0)} $. Later on, we will see that the endomorphism
$ \J:=\nabla_{\cdot }\xi + \frac12 \nabla_{J\xi}J$,
corresponding to $ \frac12 d\zeta^{(1,1)}$, defines a complex
structure on $  H$, too.

\begin{elem}\label{acs}
The endomorphisms $I$ and $K$ vanish on $\spann\{\xi, J\xi\}$
and define complex structures on $ H = \{\xi, J\xi\}^\perp$
compatible with the metric $ g$.
Moreover they satisfy $ K=I \circ J, $ 
$ 0 = I \circ J + J \circ I $ 
and for any $  X, Y \in H$ the equation
$$
(\nabla_XJ)Y
  =  
\la Y,   I X \ra \xi
  +  
\la Y,   K X \ra   J\xi .
$$
\end{elem}

We will call endomorphisms of $TM$, which are complex
structures on $ H = \{\xi, J\xi\}^\perp$ {\it transversal
complex structures}.
The last equation shows that $ \nabla J $ vanishes in
the direction of $H$. It turns out that the  same is true 
for transversal complex structures $I$ and $J$.

\begin{elem}\label{covtrans}
The transversal complex structures $  I $ and $ K $ are
parallel in the direction of the distribution $ H$, {\em i.e.}
$$
\la (\nabla_XI)Y,   Z \ra   =   
\la (\nabla_XK)Y,   Z \ra   =  
0
$$
holds for all vector fields $ X, Y, Z $ in $H$.
\end{elem}
{\sc Proof.}
First of all we compute for any vector fields $X, Y, Z$ 
the covariant derivative of $I$.
$$
\la (\nabla_X I)   Y,  Z \ra
  =  
\la (\nabla_X(\nabla_\xi I)) Y,  Z \ra
  =  
\la (\nabla^2_{X, \xi} I)Y,  Z \ra  +  
\la (\nabla_{\nabla_X\xi}I)Y,  Z \ra .
$$
If $X$ is a vector field in $ H$, then Lemma~\ref{ortho}
implies that $ \nabla_XJ $ maps $H$ to $V$ and vice versa.

Let now $X, Y, Z$ be any vector fields in $H$. Then  both summands 
in the above  formula for $\nabla_X I$ vanish, which is clear
after rewriting the first summand using formula $(5)$ of Lemma~\ref{nk} 
and the second by using  formula $(1)$ of Lemma~\ref{nk}.

The proof for the transversal complex structure $K$ is similar.

\r

Our next goal is to show that $ \frac12 d\zeta^{(1,1)} $
defines an complex structure on the orthogonal
complement of $ \spann \{\xi, J\xi\}$.

\begin{elem} \label{norm}
Let $\xi$ be a Killing vector field of length $1$ then
$$
\|d \zeta^{(1,1)}\|^2 = 8 \qquad
\|d \zeta^{(2,0)}\|^2 = 2 
$$
\end{elem}
{\sc Proof.}
We already know that the skew--symmetric endomorphism $K$ corresponding
to $-d \zeta^{(2,0)}$ is an complex structure on $\{\xi, J\xi\}^\perp$
with $K(\xi)=K(J\xi)=0$. Hence it has norm $4$ and $\|d \zeta^{(2,0)}\|^2
=\frac12 \|K\|^2=2$. Next, we have to compute the norm of 
$d \zeta^{(1,1)}$. Since $\xi $ is a Killing vector field we have 
$ \Delta \zeta = 10 \zeta$, hence
$$
\|d \zeta^{(1,1)}\|^2
 = 
\|d\zeta\|^2 - \|d \zeta^{(2,0)}\|^2
 = 
(\Delta \zeta, \zeta) - 2
 = 8
$$

\r

\begin{ecor}
The square norm of the endomorphism $\J$ corresponding to 
$ \frac12d \zeta^{(1,1)}$ is equal to $4$.
\end{ecor}

\begin{epr}
Let $(M^6, g, J)$ be a compact nearly K{\"a}hler, non--K{\"a}hler
manifold of constant type $1$ and let $\xi$ be a Killing vector field
of constant length $1$ with dual $1$--form $\zeta$. Then
$$
d^*(J\zeta) = 0 \qquad \mbox{and} \qquad \Delta(J\zeta) = 18 J\zeta .
$$
In particular, the vector field $ J\xi $ is never a 
Killing vector field.
\end{epr}
{\sc Proof.}
We start to compute the $L^2$--norm of the function $d^*(J\zeta)$:
$$
\|d^*(J\zeta)\|^2
 = 
(d^* J \zeta, d^* J \zeta)
 = 
(\Delta (J \zeta), J \zeta)
 - (d^* d (J \zeta), J \zeta)
$$
Since $\Delta=\nabla^*\nabla+\Ric  $  on $1$--forms and $\Ric=5 \Id$ we obtain
$$
\|d^*(J\zeta)\|^2
 = 
\|\nabla (J \zeta)\|^2 + 5 \|J \zeta\|^2 - \|d J \zeta\|^2
$$
To compute the norm of $\nabla (J \zeta)$ we use the formula 
$2 (\nabla X^{\flat}) = dX^{\flat} + L_Xg$ which holds for any vector
field $X$.

Note that the decomposition 
$T^*M\otimes T^*M\cong \Lambda^2(TM)\oplus\Sym^2(TM)$ 
is orthogonal. Together with Lemma~\ref{lie} this yields
\bea
\|\nabla (J \zeta)\|^2
 = 
\frac14 \left( \|d J \zeta\|^2 + \|L_{J\xi} g\|^2\right)
 = 
\frac14 \left( \|d J \zeta\|^2 + 4 \|\J\|^2\right)
\eea
From Lemma \ref{norm} follows: 
$\|\J\|^2 = \frac12 \|d \zeta^{(1,1)}\|^2 = 4$.
Using Corollary \ref{dxi}, the formula 
$d \Omega = 3 \nabla \Omega$ and the fact that 
$I=\nabla_{\xi}J$ is again  a transversal complex structure,
we find 
$$
\|d J\zeta\|^2 = 9 \|\xi \i \nabla \Omega\|^2 = 36.
$$
Combining all these computations yields $\|d^*(J\zeta)\|^2 = 0$.

Next we want to compute $\Delta J\zeta$. We start by using the Weitzenb{\"o}ck
formula on $1$--forms and the equation 
$
\nabla^*\nabla (J\tau) = (\nabla^*\nabla J)\tau + J(\nabla^*\nabla \tau) - 
2 \sum (\nabla_{e_i}J)(\nabla_{e_i}\tau)
$ for any $1$--from $\tau$.
This gives
\bea
\Delta J\zeta
&=&
\nabla^*\nabla J\zeta + \Ric(J\zeta)
 = 
(\nabla^*\nabla J)\zeta + J(\nabla^*\nabla \zeta) + 5 (J\zeta)
 - 2 \sum (\nabla_{e_i}J)(\nabla_{e_i}\zeta)\\
&=&
14 J \zeta - 2 \sum (\nabla_{e_i}J)(\nabla_{e_i}\zeta)
\eea
For the last equation we used Corollary \ref{nabla} and the assumption
that $\xi$ is 
a Killing vector field, hence $\nabla^*\nabla\zeta=\Delta
\zeta-\Ric(\zeta)=5 \zeta$.
Since $d^*(J\zeta) = 0$ we can compute $\Delta J\zeta$ by
\bea
\Delta J\zeta
&=&
d^*d(J\zeta)
 = 
\sum e_i \i \nabla_{e_i}(\xi \i d\Omega)
 = 
\sum e_i \i (\nabla_{e_i} \xi) \i d\Omega
 + e_i \i \xi \i (\nabla_{e_i} d\Omega)\qquad \\
&=&
- 3 \sum (\nabla_{e_i}\Omega)(\nabla_{e_i}\xi) + 12 J\zeta
\eea
Comparing these two equations for $\Delta J\zeta$ we get
$\sum (\nabla_{e_i}\Omega)(\nabla_{e_i}\xi) = - 2 J\zeta$,  
so finally  $\Delta J\zeta = 18 J\zeta$.

Since $\Delta X=2\Ric(X)=10X$ for every Killing vector field $X$,
$J\zeta$ cannot be Killing.

\r

Notice that if the manifold $M$ is not assumed to be compact a local
calculation still shows that $d^*(J\zeta)$ is a constant.

\begin{ecor}
If $ \xi $ is a Killing vector field of unit length, then
$$
\la d\zeta,  \Omega\ra  =  \la d\zeta^{(1,1)},  \Omega\ra  =  0 .
$$
\end{ecor}
{\sc Proof.}
Using Corollary \ref{elementary} we calculate
$$
\la d\zeta, \Omega \ra   \vol
 = 
d\zeta \wedge \ast \Omega
 = 
\frac12 d\zeta \wedge \Omega \wedge \Omega
 = 
\frac12 d(\zeta \wedge \Omega \wedge \Omega)
 = 
- d(\ast J\zeta)
 = 
\ast (d^*J\zeta)
 = 
0
$$
This proves the corollary since the decomposition 
$\Lambda^2(TM)=\Lambda^{inv}\oplus\Lambda^{anti}$
is orthogonal and $\Omega\in\Lambda^{inv}$.

\r

We are now ready to prove that the $ (1,1)$--part of $ d \zeta $ defines
a fourth complex structure on $  H =\{\xi, J\xi\}^\perp$. 
Indeed we have
\begin{epr}\label{bquadrat}
Let $\xi$ be a Killing vector field of constant length $1$. Then
$$
\J
  :=  
\nabla_{\cdot }\xi  +  \frac12 K
$$
defines a transversal complex structure on $ H $ which 
is compatible with the metric $g$.
Moreover, $  \J $ is the skew--symmetric endomorphism 
corresponding to $ \frac12d \zeta^{(1,1)} $ and it commutes with $I,\
J$ and $K$:
$$
[  \J, J  ]
  =  
[  \J,   K  ]
  =  
[  \J, I  ] =0.
$$
\end{epr}
{\sc Proof.}
Lemma~\ref{dxi20} shows that $ \J $ corresponds to $ \frac12
d\zeta^{(1,1)} $. 
Since $\la d\zeta^{(1,1)}, \Omega\ra = 0 $ and since 
$ \J $ vanishes on $\spann\{\xi, J\xi\}$,
it follows that 
$ d\zeta^{(1,1)}\in\Lambda^{(1,1)}_0(H)=\Lambda^2_-(H)$. 
Hence
$$
\J^2 = - \frac14 \|\J\|^2 \Id_H = - \Id_H.
$$
Finally, $ \J $, commutes with $I,\
J$ and $K$ since endomorphisms corresponding to self--dual and
anti--self--dual $2$--forms in 4 dimensions commute.

\r

\section{Projectable Tensors}

In this section we want to study which of the above defined tensors descend 
to the space of leaves of the integrable distribution $V = \spann
\{\xi, J\xi\}$. 
For doing this we have to compute Lie derivatives in the direction
of $\xi$ and $J\xi$. We first remark that the flow of $\xi$ preserves
both the metric $g$ and the almost complex structure $J$, thus it
also preserves $I,K,d\zeta,\J$ etc.

The situation is more complicated for $J\xi$. Since $J\xi$ is not
Killing for $g$, we have to specify, when computing the Lie
derivative of a tensor with respect to $J\xi$, whether the given
tensor is regarded as endomorphism or as bilinear form.
Here we need the following lemma, a direct consequence of the
definition 
of the Lie 
derivative and Proposition \ref{bquadrat}.

\begin{elem}\label{lie}
Let $\alpha\in\Gamma(T^*M\otimes T^*M)$ be a $(2,0)$--tensor and
$A$ be the corresponding endomorphism. Then the Lie derivatives of $A$
and $\alpha$
with respect to $J\xi$ are related by
\beq\label{li1}
(L_{J\xi} \alpha)(X, Y) = 
g((L_{J\xi}A)X, Y) + 2 g(J  \J  A(X), Y)  .
\eeq
In particular, the Lie derivative of the Riemannian metric $g$ with
respect to $J\xi$ is:
\beq\label{li2}
L_{J\xi} g=2  g(J  \J  \cdot, \cdot).
\eeq
\end{elem}

{\sc Proof.} Taking the Lie derivative in $\a(X,Y)=g(AX,Y)$ yields
$$(L_{J\xi} \alpha)(X, Y)=(L_{J\xi} g)(AX,Y)+g((L_{J\xi}A)X, Y).$$
Thus (\ref{li2}) implies (\ref{li1}). Taking $\a=g$ in (\ref{li1})
yields (\ref{li2}), so the two assertions are equivalent.

Using Proposition \ref{bquadrat} we can write:
\bea
L_{J\xi} g(X,Y)&=&g(\n_XJ\xi,Y)+g(X,\n_YJ\xi)\\
&=&\n J(X,\xi,Y)+\n J(Y,\xi,X)+g(J\n_X\xi,Y)+g(J\n_Y\xi,X)\\
&=&g(J(\J-\frac12 K)X,Y)+g(J(\J-\frac12 K)Y,X)=2g(J\J X,Y).
\eea

\r

From Lemma \ref{l35} we see that $L_{J\xi}\xi=L_{J\xi}J\xi=0$.
Thus, if $\zeta$ and $J\zeta$ denote the metric duals of $\xi$ and $J\xi$,
(\ref{li2}) shows that $L_{J\xi}\zeta=L_{J\xi}J\zeta=0$. We thus get 
\beq\label{ll1}
L_{J\xi}(d\zeta) = L_{J\xi}(d J \zeta) = 0.
\eeq

\begin{elem}\label{lieformulas}
If $\xi$ is a Killing vector field of constant length $1$, then
\beq\label{ll2}
L_{J\xi}(\Omega) = J\xi \i  d\Omega - d\zeta
 =  4\omega_K  - 2\omega_{\J} 
,\qquad
L_{J\xi}(J) = 4 K,
\eeq
\beq\label{ll3}
L_{J\xi}(\omega_K) = -4\Omega  +  4 \zeta \wedge J\zeta
 = -4\omega_J
,\qquad
L_{J\xi}(K) = -4J|_H- 2 I\J 
\eeq
\beq\label{ll4}
L_{J\xi}(\omega_{\J}) = -2\Omega  + 2 \zeta \wedge J\zeta,\qquad
L_{J\xi}(\J)=0
 \eeq
\beq\label{ll5}
L_{J\xi}(\omega_I) =  0,\qquad
L_{J\xi}(I) =  2  \J K
\eeq
where $ \omega_I,  \omega_K, \omega_{\J}$ are the $2$--forms 
corresponding to $ I,  K, \J $ and $ \omega_J $
denotes the projection of $ \Omega $ onto 
$ \Lambda^2H$, {\em i.e.} $ \omega_J = \Omega - \xi \wedge J\xi$.
\end{elem}
{\sc Proof.}
We will repeatedly use the formula 
$ L_X \alpha = X \i d \alpha + d X \i \alpha$,
which holds for any vector field $X$ and any differential form
$\alpha$.

Proposition \ref{bquadrat} shows that
\beq\label{ess}
2\omega_{\J}=d\zeta+\omega_K\eeq
hence
$$
L_{J\xi} \Omega
  =  
J\xi \i d\Omega  + d(J\xi \i \Omega)
  =  
J\xi \i d\Omega  - d\zeta=3\omega_K-d\zeta=4\omega_K-2\omega_{\J}.
$$
The second equation in (\ref{ll2}) follows by taking $A=J$ in Lemma~\ref{lie}.

Using Corollary~\ref{djxi} we get
$$L_{J\xi}(\omega_K)=
\frac13L_{J\xi} ( J\xi \i d\Omega)
  =  
\frac13J\xi \i d(J\xi \i d\Omega)
  =  
-4   J\xi \i  (J\zeta \wedge \Omega)
  =  
-4 \Omega  +  4 \zeta \wedge J\zeta.
$$
The second equation in (\ref{ll3}) 
follows directly from the first one, by taking
$A=K$ in Lemma \ref{lie}.

Using (\ref{ll1}), (\ref{ll3}) and (\ref{ess}) we obtain 
$$2L_{J\xi}\omega_{\J}=L_{J\xi}(d\zeta+\omega_K)=L_{J\xi}\omega_K=
-4 \Omega  +  4 \zeta\wedge J\zeta.$$
The second part of (\ref{ll4}) follows from Lemma \ref{lie}.

Finally, in order to prove (\ref{ll5}) we use (\ref{ll1}) twice:
$$L_{J\xi}\omega_I=\frac13 L_{J\xi}(\xi\i d\Omega)=\frac13 \xi\i
L_{J\xi}(d\Omega)=-4\xi\i(J\zeta\wedge\Omega)=0,$$
and the second part follows from Lemma~\ref{lie} again.

\r


\section{The Transversal Involution}

We define a transversal orthogonal involution $  \sigma \in \End(TM)  $ by
$$
\sigma   =   K \circ \J,
$$
{\em i.e.} we have $ \sigma^2 = \Id $ on $ H $ and $ \sigma = 0 $
on $ V=\spann\{\xi, J\xi\}$. Hence, the distribution $ H $ splits into
the $(\pm 1)$--eigenspaces of $  \sigma $ and we can define a new
metric $ g_0 $ on $ M $ as
\begin{equation}\label{g0}
g_0 = g + \frac12 g(\sigma \cdot, \cdot),
\end{equation}
{\em i.e.} we have $g_0 = g$ on $V$, $g_0 = \frac12 g $ on the
$(-1)$--eigenspace of $ \sigma $ and $ g_0 = \frac32 g $ on the
$(+1)$--eigenspace of $ \sigma$. The reason for introducing $g_0$ is
the fact that, in contrast to $g$,  this new
metric is preserved by the flow of $J\xi$ (cf. Corollary \ref{cor47} below).

\begin{elem} If $A^{\flat}$ denotes the $(2,0)$--tensor corresponding
to an endomorphism $A$, and $\alpha^{\sharp}$ denotes the endomorphism
corresponding to a $(2,0)$--tensor $\alpha$ with respect to the metric
$g$, then
$$
L_{J\xi}   \sigma^{\flat}=  - 4 (J   \J)^{\flat}
$$
\end{elem}
{\sc Proof.}
The right parts of (\ref{ll3}) and (\ref{ll4}) read
$$L_{J\xi}(K) = -4J  +  4 (\zeta\wedge\J\zeta)^{\sharp}- 2 I\J ,
\qquad
L_{J\xi}(\J)=2 (\zeta\wedge\J\zeta)^{\sharp}.$$
We clearly have $(\zeta\wedge\J\zeta)^{\sharp} \circ K=K \circ
(\zeta\wedge\J\zeta)^{\sharp}=(\zeta\wedge\J\zeta)^{\sharp} \circ
\J=\J \circ(\zeta\wedge\J\zeta)^{\sharp}=0$, 
therefore
$$
L_{J\xi}   \sigma
  =  
(L_{J\xi}  K)  \J   +   K  (L_{J\xi}  \J)
  =  
(- 4 J +4(\zeta\wedge\J\zeta)^{\sharp} - 2 J \J K)  \J
  =  
-4 J \J - 2 I.
$$
Thus Lemma~\ref{lie} gives
$$
L_{J\xi}   \sigma^{\flat}
  =  
(-4 J \J - 2 I)^{\flat}   +   2 (J \J ( K \J))^{\flat}
  =  
- 4 (J   \J)^ {\flat}.
$$
\r

\begin{ecor}\label{cor47} The metric $g_0$ is preserved by the flow of
  $J\xi$:
$$
L_{J\xi}   g_0=0.
$$
\end{ecor}
{\sc Proof.}
Direct consequence of (\ref{li2}):
$$L_{J\xi}   g_0
  =   L_{J\xi} (g + \frac12  \sigma^{\flat})=2(J\J)^{\flat}+\frac12(- 4 (J
  \J)^{\flat})=0.
$$
\r

\begin{epr}\label{lev}
For every horizontal vector fields $X,\ Y,\ Z\in H$, the Levi--Civita 
connection $ \nabla^{g_0} $ of $ g_0 $
is related to the Levi--Civita 
connection $ \nabla$ of $g$ by the formula
\beq\label{lv}
g_0(\nabla^{g_0}_X Y,   Z)
  =  
g_0(\nabla_X Y,  Z)
  +  
\frac13 
g_0((1 - \tfrac{1}{2} \sigma) 
[(\nabla_X \sigma)Y  +  (\nabla_{KX}\J)Y], Z) .
\eeq
\end{epr}
{\sc Proof.}
Since the expression $g_0(\nabla^{g_0}_X Y,   Z)-g_0(\nabla_X Y,  Z)$
is tensorial, we may suppose that $X,Y,Z$ are $\n$--parallel at some
point where the computation is performed.
The Koszul formula for $ g_0(\nabla^{g_0}_X Y,   Z) $ yields directly
\begin{equation}\label{g5}
2 g_0(\nabla^{g_0}_X Y,   Z)
  =  
2 g_0(\nabla_X Y,   Z)
  +  \frac{1}{2} 
\left[
\la (\nabla_X \sigma ) Y,  Z\ra 
  +  
\la (\nabla_Y \sigma ) X,  Z \ra 
  -  
\la (\nabla_Z \sigma ) X,   Y\ra 
\right] \ ,
\end{equation}
where $ \la \cdot, \cdot \ra $ denotes the metric $ g$.
Since $ \sigma = K \circ \J = \J \circ K $ and since -- according to
Lemma~\ref{covtrans} -- $ K $ is $\n$--parallel
in direction of $ H $,  we obtain 
$
\la (\nabla_X \sigma ) Y,  Z\ra
  =  
\la (\nabla_X \J)   K   Y,  Z\ra .
$
Now, (\ref{ess}) shows that 
$ 2\omega_{\J} - \omega_K = d\xi $ is a closed $2$--form. Hence,
$$
\la (\nabla_{X_1} \J) X_2,  X_3\ra
  +   
\la (\nabla_{X_2} \J) X_3,  X_1\ra
  +   
\la (\nabla_{X_3} \J) X_1,  X_2\ra 
  =   0
$$
for all vectors $X_i$.
Using this equation for $ X_1 = Y, X_2 = KX $ and
$ X_3 = Z $ we obtain
\bea
\firstline{1cm}{
\la (\nabla_X \sigma ) Y,  Z\ra 
  +  
\la (\nabla_Y \sigma ) X,  Z \ra 
  -  
\la (\nabla_Z \sigma ) X,   Y\ra 
  =  }
\\[1.5ex]
&=&
\la (\nabla_X \J) K Y,  Z\ra 
  +  
\la (\nabla_Y \J) K X,  Z \ra 
  -  
\la (\nabla_Z \J) K X,   Y\ra 
\\[1.5ex]
&=&
\la (\nabla_X \J) K Y    +    (\nabla_{KX} \J) Y,   Z \ra
  =  
\la (\nabla_X \sigma)  Y    +    (\nabla_{KX} \J) Y,   Z \ra
\eea
The desired formula then follows from (\ref{g0}) and (\ref{g5})
using
\begin{equation}\label{idh}
(\id_H + \frac12 \sigma)(\id_H - \frac12 \sigma ) = \frac34 \id_H.
\end{equation}

\r

We consider the space of leaves, denoted by $N$, of the integrable 
distribution $V = \spann \{\xi, J\xi\}$. The $4$--dimensional 
manifold $N$ is {\em a priori} only locally defined. It can be thought
of as the base space of  
a locally defined principal torus bundle $\mathbb{T}^2 \hookrightarrow
M \to N$. The local  
action of the torus is obtained by integrating the vector fields $\xi$ and 
$\xi^{\prime}=\frac{1}{2\sqrt{3}}J\xi$. Moreover, if one considers the
$1$--forms $\zeta$ and $\zeta'$ on $M$  
associated via the metric $g$ to the vector fields $\xi$ and
$2\sqrt{3} J\xi$ it follows that  
$\zeta(\xi)=\zeta'(\xi^{\prime})=1$ and the Lie derivatives of
$\zeta$ and $\zeta'$ in the directions  
of $\xi$ and $\xi^{\prime}$ vanish by (\ref{ll1}). 
Therefore $\zeta$ and $\zeta'$ are principal connection 
$1$--forms in the torus bundle $\mathbb{T}^2 \hookrightarrow M \to N$.

A tensor field on
$M$ projects to $N$ if and only if  it is horizontal and its Lie
derivatives with respect to $\xi$ and 
$J\xi$ both vanish. All horizontal tensors defined above have vanishing Lie
derivative with respect to $\xi$. Using (\ref{ll5}) together with 
Corollary \ref{cor47} we see that $\o_I$ and $g_0$ project down to
$N$. Moreover, $\o_I$ is compatible with $g_0$ in the sense that 
$$ \omega_I(X,Y)=\frac{2}{\sqrt{3}}g_0(I_0 X, Y),\qquad\forall\
X,Y\in H,$$ 
where $I_0$ is the $g_0$--compatible complex structure on $H$ given by 
$$
I_0=\frac{2}{\sqrt{3}}(I-\frac{1}{2}\sigma I). 
$$
This follows directly from (\ref{g0}) and (\ref{idh}). Keeping the
same notations for the projections on $N$ of projectable tensors (like
$g_0$ or $I_0$) we now prove

\begin{ath} \label{mt}
$(N^4,g_0,I_0)$ is a K{\"a}hler manifold.
\end{ath}
{\sc Proof.}
In order to simplify notations we will denote by $\t\nabla$ and 
$\t\nabla^{g_0}$ the
partial connections on the distribution $H$ given by the
$H$--projections of the Levi--Civita  
connections $\nabla$ and $\nabla^{g_0}$. 

Then Proposition \ref{lev} reads 
\beq\label{lle1}
\t\nabla^{g_0}_X = \t\nabla_X +  
\frac13 (\id_H-\frac12\sigma)(\t\n_X\sigma+\t\n_{KX}\J).
\eeq

We have to check that $\t\nabla^{g_0}_XI_0=0$
for all $X$ in $H$. We first notice the tautological relation
$\t\nabla_X\id_H =0$.
From Lemmas \ref{acs} and \ref{covtrans} we have 
$\t\nabla_XI=\t\nabla_XJ=\t\nabla_XK=0$ 
for all $X$ in $H$. Moreover, the fact that
$I,\ J$ and $K$ commute with $\J$ and the relation $\J^2=\id_H$ easily
show that $\t\n_X\J$ commutes with $I,\ J,\ K$ and anti--commutes with $\J$
and $\s$. Consequently, $\t\n_X\s\ (=K\t\n_X\J)$ commutes with $K$ and 
anti--commutes with $I,\ J,\ \J$ and $\s$ for all $X\in H$.

We thus get
\beq\label{ji1}
 \t\nabla_X I_0=\t\nabla_X\frac{2}{\sqrt{3}}(I-\frac{1}{2}\sigma I)
=-\frac{1}{\sqrt{3}}(\t\nabla_X\s)I.
\eeq
On the other hand, the commutation relations above show immediately
that the endomorphism $I_0$ commutes with
$(\id_H-\frac12\sigma)\t\n_{KX}\J$ and 
anti--commutes with $(\id_H-\frac12\sigma)\t\n_X\sigma$. Thus 
the endomorphism $\frac13
(\id_H-\frac12\sigma)(\t\n_X\sigma+\t\n_{KX}\J)$
 acts on $I_0$ by 
\bea\frac13
(\id_H-\frac12\sigma)(\t\n_X\sigma+\t\n_{KX}\J)(I_0)&=&
2\frac13
(\id_H-\frac12\sigma)(\t\n_X\sigma)I_0\\
&=&\frac{4}{3\sqrt3}(\id_H-\frac12\sigma)
\t\n_X\sigma(\id_H-\frac12\sigma)I\\
&=&\frac{4}{3\sqrt3}(\id_H-\frac12\sigma)(\id_H+\frac12\sigma)
(\t\n_X\sigma)I\\
&=&\frac{1}{\sqrt{3}}(\t\nabla_X\s)I. 
\eea

This, together with (\ref{lle1}) and (\ref{ji1}), shows that 
$\t\nabla^{g_0}_XI_0=0$.

\r

We will now look closer at the structure of the metric $g$. Since 
$$g=g_0+\frac12 \omega_K(\J\cdot, \cdot),$$
the geometry of $N$, together with the form $\omega_K$ and the almost
complex structure $\J$ determine
completely the nearly K{\"a}hler metric $g$. But the discussion below will 
show that $\omega_K$ depends also in an explicit way on the geometry of
the K{\"a}hler surface $(N^4,g_0, I_0)$.

If $\alpha$ is a $2$--form on $H$ we shall denote by $\alpha^{\prime}$,
resp. $\alpha^{\prime \prime}$ the invariant, resp. anti--invariant parts of
$\alpha$ with respect to the almost complex structure $I_0$. An easy algebraic
computation shows that $\omega_J$ is $I_0$--anti--invariant whilst 
\begin{equation}\label{52}
 \omega_K^{\prime}=-\frac{1}{3}(\omega_K-2\omega_{\J}), \
\mbox{and} \ \omega_K^{\prime \prime}=\frac{2}{3}(2
\omega_K-\omega_{\J}).
\end{equation}
Consider now the complex
valued $2$-form of $H$ given by  
\begin{equation}
 \Psi=\sqrt{3}\omega_K^{\prime \prime}+2i\omega_J. 
\end{equation}
It appears then from Lemma \ref{lieformulas} and (\ref{52}) that 
\begin{equation}\label{60}
 L_{\xi^{\prime}}\Psi=i\Psi. 
\end{equation}

Thus $\Psi$ is not projectable on $N$, but it can be interpreted as a
$\L$--valued $2$--form on $N$, where $\L$ is the complex line bundle
over $N$ associated 
to the (locally defined) principal $S^1$--bundle 
$$M/\{\xi\}\to N:=M/\{\xi,\xi'\}$$
with connection form $\zeta'$.

Corollary \ref{dxi}, together with
(\ref{ik}) implies that the curvature form of $\L$ equals
\beq\label{dz}
d \zeta'=-6\sqrt{3} \omega_I=-12g_0(I_0
\cdot, \cdot).
\eeq
Notice that, since the curvature form of $\L$ is of type $(1,1)$, the
Koszul--Malgrange theorem implies that $\L$ is holomorphic.

The following proposition computes the Ricci curvature of the K{\"a}hler
surface $(N^4,g,I_0)$ by identifying the line bundle $\L$ with the
anti--canonical bundle of $(N,I_0)$.

\begin{epr}\label{kahl}
$(N^4,g,I_0)$ is a K{\"a}hler-Einstein surface with Einstein constant equal to
$12$. Moreover, $\L$ is isomorphic to the anti--canonical line bundle
$\K$ of $(N^4,g,I_0)$.
 \end{epr} 
{\sc Proof.}
We first compute $\omega_J(I_0 \cdot, \cdot)=-\frac{\sqrt{3}}{2}
\omega_K^{\prime \prime}$ and $(\omega_K^{\prime \prime})(I_0 \cdot, \cdot)=
\frac{2}{\sqrt{3}} \omega_J$. These lead to 
$$ \Psi(I_0 \cdot, \cdot)=-i \Psi $$
in other words $\Psi$ belongs to $\Lambda^{0,2}_{I_0}(H, \mathbb{C})$. 
We already noticed that by (\ref{60}), $\Psi$ defines a
section of the holomorphic line bundle
\begin{equation}\label{59}
 \Lambda^{0,2}_{I_0}(N)\otimes\L=\K^{-1} \otimes\L. 
\end{equation}
Since $\Psi$ is non--vanishing, this section induces an isomorphism
$\Psi: \K\to\L$. We
now show that $\Psi$ is in fact $\t\n^{g_0}$--{\it{parallel}}. 

Notice first that $K$ commutes with
$(\id_H-\frac12\sigma)(\t\n_X\sigma+\t\n_{KX}\J)$ (it actually
commutes with each term of this endomorphism), and $\t\n K=0$ by Lemma
\ref{covtrans}. Thus (\ref{lle1}) shows that 
\beq\label{t1}\t\n^{g_0}K=0.
\eeq 
Furthermore, using the relation 
$$g(\cdot,\cdot)=\frac43
g_0((1-\frac{\sigma}2)\cdot,\cdot),$$
 $\Psi$ can be expressed as
\beq\label{t2}\Psi=\frac{4}{\sqrt 3}g_0((K-iI_0 K)\cdot,\cdot).\eeq
Since $\t\n^{g_0} g_0=0$ and $\t\n^{g_0}I_0=0$ (by Theorem \ref{mt}),
(\ref{t1}) and (\ref{t2}) show
that $\Psi$ is $\t\n^{g_0}$--parallel. 

Hence the Ricci form of $(N^{4},g_0,I_0)$ is opposite to the 
curvature form of $\L$. From (\ref{dz}) we obtain 
$Ric_{g_0}=12g_0$, thus finishing the proof.   

\r

\begin{epr} The almost complex structure $\J$ on $H$ is projectable and
  defines an almost K{\"a}hler structure on $(N,g_0)$ commuting with $I_0$.
\end{epr}

{\sc Proof.}
Lemma \ref{lieformulas} shows that $\J$ is projectable onto $N$. Let
us denote the
associated $2$--form with respect to $g_0$ by $\o^0_{\J}$. 
Identifying forms and endomorphisms via the metric $g$ we can write
$$\o^0_{\J}:=g_0(\J\cdot,\cdot)
=(1+\frac12\sigma)\J=(1+\frac12 K\J)\J=\J-\frac12 K=\frac12d\zeta.$$
This shows that $\o^0_{\J}$ is closed, so the projection of 
$(g_0,\J)$ onto $N$ defines an almost
K{\"a}hler structure.

\r

Together with Proposition \ref{kahl}, we see 
that the locally defined manifold $N$ carries
a K{\"a}hler structure $(g_0,I_0)$ and an almost K{\"a}hler
structure $(g_0,\J)$, both obtained by projection from $M$. 
Moreover $g_0$ is Einstein with positive scalar curvature. If $N$ were
compact, we could directly apply Sekigawa's proof of the Goldberg
conjecture in the positive curvature case in order to conclude that
$(g_0,\J)$ is K{\"a}hler. As we have no information on
the global geometry of $N$, we use the following idea.
On any almost K{\"a}hler Einstein manifold, a Weitzenb{\"o}ck--type formula was 
obtained in \cite{adm}, which in the compact case shows by integration 
that the manifold is actually K{\"a}hler provided the Einstein constant
is non--negative. In the present situation, we simply interpret on $M$
the corresponding formula on $N$, and after integration over $M$ we prove a
pointwise statement which down back on $N$ just gives the
integrability of the almost K{\"a}hler structure.

The following result is a particular case (for Einstein metrics) of 
Proposition 2.1 of \cite{adm}:

\begin{epr} \label{prop2} For any almost K{\"a}hler Einstein manifold
$(N^{2n},g_0,J,\O)$ with covariant derivative denoted by $\n$ and
curvature tensor $R$, the following pointwise relation holds:

\begin{equation}\label{la1}
\Delta^{N} s^*-
8 \delta^{N} (\langle \rho^* , \n_{\cdot}  \O\rangle) = -
8|R''|^2  - |\n^*\n \O |^2 - |\phi|^2
-\frac{s}{2n}|\n\O|^2,
\end{equation}
where $s$ and $s^*$ are respectively the scalar and $*$--scalar
curvature, $\rho ^*:=R(\O)$ is the $*$--Ricci form, 
$\phi(X,Y) = \langle \n_{JX} \O,\n_Y \O \rangle $,
and $R''$ denotes the projection of the curvature tensor on the
space of endomorphisms of $[\Lambda ^{2,0}N]$ anti--commuting with $J$. 
\end{epr}

We apply this formula to the (locally defined) 
almost K{\"a}hler Einstein manifold
$(N,g_0,\J)$ with Levi--Civita covariant derivative denoted $\n^0$ and
almost K{\"a}hler form $\hat\O$ and obtain
\begin{equation}\label{fdn}F+\delta ^N \alpha =0,
\end{equation}
where 
$$F:=8|R''|^2+|(\n^0)^*\n^0 \hat\O
|^2+|\phi|^2+\frac{s}{4}|\n^0\hat \O|^2$$
is a non--negative function on $N$ and 
$$\alpha:=ds^*-8g_0(\rho^*,
\n^0_{\cdot}\hat \O)$$
is a 1--form, both $\alpha$ and $F$ depending in an explicit way on
the geometric data $(g_0,\J)$. Since the Riemannian submersion 
$\pi:(M,g_0)\to N$ has minimal (actually totally geodesic) fibers, the 
codifferentials on $M$ and $N$ are related by $\delta^M(\pi ^*\a)=\pi ^*
\delta^N\a$ for every 1--form $\a$ on $N$. Thus (\ref{fdn}) becomes
\begin{equation}\label{fdn1}\pi ^*F+\delta ^M (\pi ^*\alpha) =0.
\end{equation}
Notice that the function $\pi ^*F$ and the $1$--form $\pi ^*\alpha$ are
well--defined {\em global} objects on $M$, even though $F$, $\a$ and
the manifold $N$ itself 
are just local. This follows from the fact that $F$ and
$\a$ only depend on the geometry of $N$, so $\pi ^*F$ and $\pi
^*\alpha$ can be explicitly defined in terms of $g_0$ and $\J$ on $M$.

When $M$ is {\it{compact}}, since $\pi ^*F $ is non--negative,
(\ref{fdn1}) yields, after 
integration over $M$, that $\pi ^*F=0$. Thus $F=0$ on $N$ and this
shows, in particular, that $\phi=0$, so $\J$ is parallel on $N$.

\section{Proof of Theorem \ref{main}}

By the discussion above, when $M$ is compact, $\J$ is parallel on $N$
with respect to the Levi--Civita 
connection of the metric $g_0$, so $\J$ is
$\t\n^{g_0}$--parallel on $H$.

\begin{elem} \label{l61} The involution $\s$ is $\t\n$--parallel.
\end{elem}
{\sc Proof.} Since $\s=\J K$, (\ref{t1}) shows that 
$\tilde{\nabla}^{g_0}\s=0$. Using (\ref{lle1}) and the fact that $\s$
anti--commutes with  $\tilde{\nabla}_X\sigma$ and $\tilde{\nabla}_X\J$
for every $X\in H$, we obtain
$$ \tilde{\nabla}_X\s+\frac{2}{3} 
(id_H-\frac{1}{2}\sigma)(\tilde{\nabla}_X\sigma+\tilde{\nabla}_{KX}\J)
\s=0$$ 
for all $X$ in $H$. Since $I$ commutes with $\tilde{\nabla}_X\J$ and 
anti--commutes with $\s$ and $\tilde{\nabla}_X\sigma$, the
$I$--invariant part of the above equation reads
\beq\label{pp}\frac23
(\tilde{\nabla}_X\sigma)\s+\frac13\tilde{\nabla}_{KX}\J=0.\eeq 
But $\s=\J K$ and $\t\n K=0$, so from (\ref{pp}) we get
$$2(\t\n_X\J)\J=\tilde{\nabla}_{KX}\J.$$
Replacing $X$ by $KX$ and applying this formula twice yields
$$\tilde{\nabla}_{X}\J=-2(\t\n_{KX}\J)\J=4\t\n_X\J,$$
thus proving the lemma.

\r

We now recall that the first canonical Hermitian connection of the
NK structure $(g,J)$ is given by 
$$ \overline{\nabla}_U=\nabla_U+\frac{1}{2}(\nabla_U J)J $$
whenever $U$ is a vector field on $M$. We will show that $(M^6,g)$ is
a homogeneous space actually by showing that 
$\overline{\nabla}$ is a Ambrose--Singer connection, that is
$\overline{\nabla} \bar T=0$ and $\overline{\nabla} \bar{R}=0$, where
$\bar{T}$  and $\bar{R}$ denote the torsion and
curvature tensor of the canonical connection $\overline{\nabla}$. 

Let $H_{\pm}$ be the eigen--distributions of the involution
$\sigma$ on $H$, corresponding to the eigenvalues $\pm 1$. We define the new
distributions
$$ E=<\xi> \oplus H_{+} \ \mbox{and} \ F=<J\xi> \oplus H_{-}.$$
Obviously, we have a $g$--orthogonal splitting $TM=E \oplus F$, with
$F=JE$. 

\begin{elem} The splitting $TM=E \oplus F$ is parallel with respect to
  the first canonical connection. 
\end{elem}
{\sc Proof. }
For every tangent vector $U$ on $M$ we can write
$$ \overline{\nabla}_U\xi=\nabla_U \xi+\frac{1}{2}(\nabla_UJ)J\xi
=\J U-\frac{1}{2}KU+\frac{1}{2}JIU= (\sigma+1)\J U, $$
showing that $\overline{\nabla}_U\xi$ belongs to $E$ (actually to
$H_+$) for all $U$ in $TM$. 

Let now $Y_{+}$ be a local section of $H_{+}$. We have to consider three
cases. First, 
$$\overline{\nabla}_{\xi} Y_{+}=\n_{\xi}Y_++\frac12 (\n_{\xi}J)JY_+
=L_{\xi}Y_++\n_{Y_+}\xi+\frac12 IJY_+=L_{\xi}Y_++\J Y_+$$
belongs to $H_+$ since $L_\xi$ and $\J$ both preserve $H_+$. Next, 
if $X$ belongs to $H $ then: 
$$ 
\overline{\nabla}_{X} Y_{+}=\tilde{\nabla}_XY_{+}+\la
\overline{\nabla}_XY_{+},\xi\ra \xi+\la 
\bar{\nabla}_XY_{+},J\xi\ra J\xi $$
But $< \overline{\nabla}_XY_{+},J\xi> =<JY_{+},
\bar{\nabla}_X\xi>=0$ by the above discussion 
and the fact that $JY_{+}$ is in $H_{-}$, and
$\t{\nabla}_XY_{+}$ is an element of $H_+$ by Lemma \ref{l61}. Thus 
$\overline{\nabla}_{X} Y_{+}$ belongs to $E$.

The third case to consider is 
\bea\overline{\nabla}_{J\xi} Y_{+}&=&{\nabla}_{J\xi} Y_{+}+\frac12
(\n_{J\xi}J)JY_+ =L_{J\xi} Y_{+}+\n_{Y_+}J\xi-\frac12\n_{\xi}Y_+\\
&=&L_{J\xi} Y_{+}+(\n_{Y_+}J)\xi+J\n_{Y_+}\xi-\frac12 IY_+\\
&=&L_{J\xi} Y_{+}-IY_++J(\J Y_+-\frac12 KY_+)-\frac12 IY_+\\
&=&L_{J\xi} Y_{+}-2IY_++J\J Y_+.
\eea
On the other hand
$$L_{J\xi}Y_+=L_{J\xi}\s Y_+=\s L_{J\xi}Y_++(L_{J\xi}\s) Y_+,$$
so the $H_-$--projection of $L_{J\xi}Y_+$ is
$$\pi_{H_-}L_{J\xi}Y_+=\frac{1-\s}{2}L_{J\xi}Y_+=\frac12 (L_{J\xi}\s)
Y_+.$$
Using (\ref{ll3}) and (\ref{ll4}) and the previous calculation we get 
\bea\pi_{H_-}(\overline{\nabla}_{J\xi} Y_{+})
&=&\pi_{H_-}(L_{J\xi} Y_{+}-2IY_++J\J Y_+)=\pi_{H_-}(\frac12 (L_{J\xi}\s)
Y_+-2IY_++J\J Y_+)\\
&=&\pi_{H_-}((I-2J\J-2I+J\J)Y_+)=-\pi_{H_-}((1+\s)(Y_+))=0.
\eea
Thus $E$ is $\overline{\nabla}$--parallel, and since $F=JE$ and
$\overline{\nabla}J=0$ by definition, we see that
$F$ is $\overline{\nabla}$--parallel, too.

\r

Therefore the canonical Hermitian connection of $(M^6,g,J)$ has
reduced holonomy, more precisely complex 
irreducible but real reducible. Using Corollary 3.1, page 487 of
\cite{Nagy} we obtain that $\overline{\nabla} \bar{R}=0$. 
Moreover, the condition $\overline{\nabla} \bar{T}=0$ is always
satisfied on a NK manifold (see \cite{bm}, lemma 2.4 for instance). The
Ambrose--Singer theorem shows that if $M$ is simply 
connected, then it is a homogeneous space. To conclude
that $(M,g,J)$ is actually $S^3 \times S^3$ we use the fact that the
only homogeneous NK manifolds 
are $S^6, S^3 \times S^3, \mathbb{C}P^3, F(1,2)$ (see
\cite{Butruille}) and 
among these spaces only $S^3 \times S^3$ has vanishing Euler
characteristic. If $M$ is not simply connected, one applies the
argument above to the universal cover of $M$ which is compact and finite by
Myers' theorem. The proof of Theorem \ref{main} is now complete.

\r

\section{The Inverse Construction} 

The construction of the (local) torus bundle $M^6\to N^4$ described in
the previous sections gives rise to the following Ansatz for constructing local
NK metrics. 

Let $(N^4,g_0,I_0)$ be a (not necessarily complete) 
K{\"a}hler surface with $Ric=12g_0$ and assume
that $g_0$ carries a compatible almost--K{\"a}hler
structure $\J$ which commutes with $I_0$. Let $\L \to N$ be
the anti--canonical line bundle of $(N^4,g_0,I_0)$
and let $ \pi_1 : M_1 \to N$ be the associated principal circle
bundle. Fix 
a principal connection form $\theta$ in $M_1$ with curvature
$-12\omega_{(g_0,I_0)}$. 
Let  $H$ be the horizontal distribution of this connection and let  $\Phi$ in
$\Lambda^{0,2}_{I_0}(H, \mathbb{C})$ be the "tautological" 
$2$--form obtained by the lift of the identity 
map $1_{\L^{-1}} : \L^{-1} \to \L^{-1}$. 

Give $M_1$ the Riemannian metric   
\begin{equation}
g_1=\theta\otimes\theta+\frac{2}{3}
\pi_1^{\star}g_0-\frac{1}{2\sqrt{3}} (Re \Phi)(\J \cdot,
\cdot). 
\end{equation}
Let now $M$ denote the principal $S^1$--bundle $ \pi : M \to M_1 $
with first Chern class represented by the closed $2$--form
$\Omega=2\pi_1^{\star}g_0(\J \cdot, \cdot)$. Since we work locally
we do not have to worry about integrability matters. Let $\mu$ be a connection
$1$--form in $M$ and give $M$ the Riemannian metric  
$$g=\mu^2+\pi^{\star}g_1.$$ 
We consider on $M$ the $2$--form  
\begin{equation}
\omega=\frac{1}{2\sqrt{3}} \mu \wedge
\pi^{\star}\theta+\frac{1}{2}\pi^{\star}(Im\Phi). 
\end{equation}

By a careful inspection of the discussion in the previous sections we
obtain:    
\begin{epr} $(M^6,g, \omega)$ is a nearly K{\"a}hler manifold of constant
type equal to $1$.
Moreover, the vector field dual to $\mu$ is a unit Killing
vector field.   
\end{epr}  

Notice that the only {\em compact} K{\"a}hler--Einstein surface
$(N^4,g_0,I_0)$ with $Ric=12g_0$ possessing an almost K{\"a}hler
structure commuting with $I_0$ is the product of two spheres of radius
$\frac1{2\sqrt3}$ (see \cite{adm}), which corresponds, by the above
procedure, to the nearly K{\"a}hler structure on $S^3\times S^3$. Thus
the new NK metrics provided by our Ansatz cannot be compact, which is
concordant with Theorem \ref{main}.

\labelsep .5cm

\end{document}